\documentclass{amsart}
\usepackage[mathscr]{eucal}
\usepackage{graphicx}
\usepackage{amscd}
\usepackage{amsfonts}
\usepackage{amsmath}
\usepackage{amsthm}
\usepackage{amssymb}
\usepackage{latexsym}
\newfont{\msam}{msam10}

\newtheorem{theorem}[]{Theorem}

\newtheorem{lemma}[]{Lemma}
\theoremstyle{definition}                               

\newtheorem{remark}[]{Remark}

\def\remark{\noindent\textbf{Remark.}} 
\newcommand{\la}{\label}
\newcommand{\into}{\,\,\hookrightarrow\,\,}
\def\C{\mathbb{C}}
\def\A{\mathbb{A}}
\def\Z{\mathbb{Z}}
\def\R{\mathbb{R}}
\def\bs#1{\boldsymbol{#1}}
\def\ms#1{\mathcal{#1}}

\def\Gm#1{\mbox{\sf GrMod}(#1)}

\def\Gma{\mbox{\sf GrMod}(\bs{A})}

\def\P{\mathbb{P}_{q}^{2}}

\def\coh{\mbox{\sf Coh}}

\def\Hom\sub#1#2{\underline{\mbox{\rm Hom}}_{\bs{A}}(#1,#2)}
\def\Ext\sub#1#2{\underline{\mbox{\rm Ext}}^{i}_{\bs{A}}(#1,#2)}
\numberwithin{equation}{section}
 \begin{document}
\title{A Remark on Letzter-Makar-Limanov Invariants}
%
%
	\author{Yuri Berest}
\address{Department of Mathematics, Cornell University, 
Ithaca, NY 14853-4201, USA}
\email{berest@math.cornell.edu}
\thanks{Research supported by the NSF grant DMS 00-71792 
and an A.~P.~Sloan Research Fellowship.}
\maketitle
\section{Introduction}

Let $ X $ be an irreducible affine algebraic curve over $ \C $ with 
normalization $ \tilde X \,$, and let $ D(X) $ be the ring of (global) 
differential operators on $ X \,$. Due to general results 
of Smith and Stafford (see \cite{SS}), it is known that $ D(X) $ is Morita 
equivalent to $ D(\tilde X) $ if (and only if) the normalization map 
$\, \pi: \tilde X \to X \,$ is injective. In the special case when 
$ X $ is rational and $ \tilde X $ is isomorphic to the affine line 
$ \A^{1} $, there exist non-isomorphic curves with isomorphic
rings of differential operators (the first examples of this kind 
were found in \cite{L}). To distinguish non-isomorphic rings
$ D(X) $ in the Morita equivalence class of $ D(\A^{1}) \,$
G.~Letzter and L.~Makar-Limanov \cite{LM} introduced a certain 
ring-theoretic invariant which they called `$\mbox{codim}\,D(X)$'.
Subsequently, a complete classification of rings of differential 
operators on affine curves was given (see \cite{K}, \cite{BW1}, 
and \cite{BW4} for a detailed exposition). This showed that for the 
above class of curves the algebras $ D(X) $ were determined (up to 
isomorphism) by a single non-negative integer $ n\,$. 
The relation between this new invariant $ n $ and the 
Letzter-Makar-Limanov one turned out to be very simple (cf. Remark 
in \cite{BW1}): $\, \mbox{codim}\, D(X) = 2n \,$. However, since 
the invariants in hand seemed to be of different nature, this 
relation looked somewhat mysterious: it did not follow {\it a priori}\, 
from the definitions, but rather was the result of comparing the 
theorems of \cite{K} and \cite{BW1} with some explicit calculations 
in \cite{LM}.

In this short note we will give a natural interpretation 
of the Letzter-Makar-Limanov (LM) invariants in the spirit of 
noncommutative projective geometry of M.~Artin et al. 
(see \cite{A}, \cite{AZ}, and \cite{SV} for a general overview
of the subject). As a result, we will prove the above relation 
in a very simple, purely homological fashion. We will introduce
and work with LM-invariants in a slightly more general setting 
(than that of \cite{LM}) and will mostly use the notation of 
\cite{BW3}. Apart from  \cite{AZ}, the latter reference also 
contains a review of all the definitions and results from 
noncommutative geometry needed for the present paper.

In the end, I would like to thank G.~Wilson for his suggestions 
and encouragement and the organizers of 
the ICRA X  for their kind invitation and support.

\section{Geometric Interpretation of LM-invariants}

Let $ A := A_{1}(\C) $ be the first complex Weyl algebra; 
$\, Q := \mbox{Frac}(A) \,$, 
its field of fractions; $\, M \,$, a finitely generated rank $1$ 
torsion-free right 
module over $ A\,$; $\, D := \mbox{End}_{A}(M)\,$, the 
endomorphism ring of $ M\,$. Since $ A $ is a simple 
Noetherian hereditary algebra, $ M $ is automatically a 
progenerator in the category $ \mbox{\sf Mod}(A) $ of 
right $A$-modules, and the ring $ D $ is Morita 
equivalent to $ A \,$.

As usual, we may (and will) identify $ M $ with an ideal of 
$ A $ (possibly fractional), and $ D $ with the subring 
of $ Q \,$:
\begin{equation}
\la{1}
D = \{\, q \in Q\ :\ q M \subseteq M\,\}\ .
\end{equation}

Let $ \bs{w}= (w_1, w_2) $ be a nonzero non-negative vector in 
$ \R^{2} $, and let $ \{Q_{\bullet}(\bs{w})\} $ denote the standard 
Dixmier filtration on $\, Q\,$ of weight $\,\bs{w}\,$ (see \cite{D}). 
We will equip all the subspaces of $ Q $ with induced filtrations. 
In particular, we set
\begin{equation}
\la{11}
A_{k}(\bs{w}) := A \cap Q_{k}(\bs{w}) \quad \mbox{and} \quad  
D_{k}(\bs{w}) := D \cap Q_{k}(\bs{w})\ ,
\end{equation}
and write $ \bs{GQ}(\bs{w}) := \bigoplus_{k \in \Z} 
Q_k/Q_{k-1} \,$, $\, \bs{GA}(\bs{w}) := \bigoplus_{k \in \Z} 
A_k/A_{k-1} \,$ and $\, \bs{GD}(\bs{w}) := \bigoplus_{k \in \Z}  
D_k/D_{k-1} $ for the associated graded rings. The natural 
inclusions $\, A_k \subset Q_k \,$, $\, D_k \subset Q_k \,$
induce then the embeddings of graded algebras $\, \bs{GA}(\bs{w}) 
\into  \bs{GQ}(\bs{w}) \,$, $\, \bs{GD}(\bs{w}) \into  \bs{GQ}(\bs{w})
 \,$. We will identify $\, \bs{GA}(\bs{w}) \,$ and $\, 
\bs{GD}(\bs{w}) \,$ with their images in $\,\bs{GQ}(\bs{w})\,$
and note that $\, \bs{GD}(\bs{w}) \,$ does not depend on the choice
of representation of $ M $ as an ideal in $ Q\,$.

The following theorem is essentially a restatement (perhaps, 
in a slightly more general form) of one of the main results of 
\cite{LM} (see {\it loc. cit.}, Proposition~2.4, and 
\cite{SS}, Sections 3.11, 3.12).
\begin{theorem}
\la{LM} 
For each nonzero $ \bs{w} \in \R^{2}_{+} \,$, we have 
$\, \bs{GD}(\bs{w}) 
\subseteq \bs{GA}(\bs{w}) \,$ in $\, \bs{GQ}(\bs{w})\,$; 
the quotient space $\, \bs{LM}(\bs{w}) := 
\bs{GA}(\bs{w})/\bs{GD}(\bs{w}) \,$ is finite-dimensional, 
its dimension being independent of $ \bs{w}\,$.
\end{theorem}
In what follows we will denote the number $\,\dim_{\mathbb{C}} \bs{LM}\,$ 
by $\, p_{D} \,$ and call it the {\it LM-invariant} of the ring $ D\,$; 
since $ p_{D} $ is independent of $ \bs{w} \,$, we will also assume 
$\,\bs{w} = (1,1)\,$ for simplicity of further considerations.

Let  $\, \bs{A} := \bigoplus_{k \in \Z} A_k \,$ and 
$\, \bs{D} := \bigoplus_{k \in \Z} D_k\,$ be the Rees algebras 
(homogenizations) of $ A $ and $ D $
with respect to the filtrations $ \{A_{\bullet}\} $ and 
$ \{D_{\bullet}\} \,$. Then $\, \bs{A} \, $ is isomorphic 
to the algebra of `noncommutative polynomials' $\,\C[x,y,z]\,$
with generators $\, x, y, z \,$, all having degree 
$ 1\,$, $ z $ commuting with $ x $ and $ y \,$, and $\, xy-yx = z^2\,$.
On the other hand, $\, \bs{GD} \,$ and $\, \bs{D} \,$ are both graded 
connected Noetherian (and hence, locally finite) algebras. 
This can be shown easily first for $\, \bs{GD} \,$ by
adapting the techniques of \cite{SS} (specifically, the proof of 
Theorem~3.12, {\it loc. cit.}), and then for $\, \bs{D} \,$
by a standard lifting argument (see, for example, \cite{Bj}, 
Appendix~III, Proposition~1.29).

Now, following \cite{AZ}, we will think of $ \bs{A} $ and $ \bs{D} $
as the homogeneous coordinate rings of {\it noncommutative}
projective schemes $\, X_{A} := \mbox{\sf Proj}(\bs{A}) 
$ and $ X_{D} := \mbox{\sf Proj}(\bs{D}) \,$.
By definition, these are the Serre quotients 
$ \coh(X_{A}) :=  \Gma/\mbox{\sf Tors}(\bs{A}) $ and
$ \coh(X_{D}) :=  \Gm{\bs{D}}/\mbox{\sf Tors}(\bs{D}) $
of the categories of (finitely generated graded right) modules 
modulo the torsion subcategories (consisting of finite-dimensional 
modules) taken with distinguished objects 
$ \ms{O}_{X_{A}} := \widetilde{\bs{A}} \in \coh(X_{A}) $ and 
$ \ms{O}_{X_{D}} := \widetilde{\bs{D}} \in \coh(X_{D})\,$.
By analogy with the commutative case, the objects of
$ \coh(X_{A}) $ and $ \coh(X_{D}) $ are referred to as
{\it coherent sheaves}, with $ \ms{O}_{X_{A}} $ and 
$ \ms{O}_{X_{D}} $ playing the role of `structure sheaves'
on $ X_A $ and $ X_D $ respectively. Note that $ X_{A} $ 
is the quantum projective plane $ \P $ defined in 
\cite{LeB} and \cite{BW3}.

The following observation is a simple consequence of Theorem~\ref{LM}.
\begin{lemma}
\la{L1} 
The Hilbert polynomial of the graded algebra $ \bs{D} $
is given by the formula
\begin{equation}
\la{3}
P_{\bs{D}}(k) = \frac{1}{2}(k+1)(k+2) - p_{D}\ ,
\end{equation}
so that $\, p_{D} = 1 - P_{\bs{D}}(0) \,$. By analogy with 
the commutative case, we may thus interpret $\, p_{D} \,$
as (minus) the {\rm arithmetic genus} of the `projective surface' 
$ X_{D} \,$.
\end{lemma}
\begin{proof}
Being finite, the graded vector space $ \bs{LM} $ is bounded from above.
Therefore, by Theorem~\ref{LM}, choosing $ k \gg 0 \,$ we can write
the exact sequence of {\it finite-dimensional} graded vector spaces:
\begin{equation}
\la{5}
0 \to \bs{GD}_{\leq k} \to \bs{GA}_{\leq k} \to \bs{LM} \to 0\ .
\end{equation}
Whence we have
\begin{eqnarray}
\dim \,\bs{LM} & = & \dim \,\bs{GA}_{\leq k} - 
 \dim \,\bs{GD}_{\leq k}  \nonumber \\
& = &  \dim \bigoplus\limits_{i=-\infty}^{k}
 A_{i}/A_{i-1} - \dim
\bigoplus\limits_{i=-\infty}^{k} D_{i}/D_{i-1} \nonumber \\
& = & \sum\limits_{i=-\infty}^{k} (\dim A_{i} - 
\dim A_{i-1}) - 
\sum\limits_{i=-\infty}^{k} (\dim D_{i} - 
\dim D_{i-1}) \nonumber \\*[1ex]
& = & \dim A_{k} - \dim D_{k}\quad \mbox{ for all}\ k \gg 0\ .\nonumber
\end{eqnarray}
By definition, $\, P_{\bs{D}}(k) = \dim D_{k} \,$ for $ k \gg 0 \,$,
and $\, \dim A_{k} = (k+1)(k+2)/2 \,$ for all $ k \geq 0\,$. Combined
together these give formula (\ref{3}), and therefore finishes the proof
of the lemma.
\end{proof}

\remark\  The category of affine schemes is dual (that is, 
anti-equivalent) to the category of commutative rings and ring 
homomorphisms. In particular, two schemes are isomorphic if and only 
if their coordinate rings are isomorphic. Passing to noncommutative 
domain it seems more appropriate to associate `a noncommutative affine 
scheme' not to an isomorphism class but to a Morita equivalence class 
of a given noncommutative ring (see, for example, \cite{Sm} and 
\cite{SV}). 
In other words, two  Morita equivalent rings should be thought of as 
coordinate rings of {\it isomorphic} noncommutative spaces. 
From this perspective the result of Lemma~\ref{L1} appears to be 
natural. Indeed, if we define the quantum plane $\, \mathbb{A}_{q}^{2}\,$ 
as an `affine space' associated with the Morita class of $ A \,$,  
the projective schemes $ X_{D} $ and $ X_{A} $ can be regarded as 
different `compactifications' of $\, \mathbb{A}_{q}^{2}\,$.
In that case {\it any} invariant that may distinguish $ A $ and 
$ D $ up to isomorphism should be an object of projective (rather 
than affine) geometry.

\section{The Relation $\, p_{D} = 2 n\,$}

First, we recall the definition of the number $ n \,$.
As shown in \cite{BW3}, Section~4, every ideal (class of) $ M $ 
admits a unique extension to a coherent rank $ 1 $ torsion-free
sheaf $\, \ms{M} \,$ on $ \P $ trivial over the line at infinity. 
This extension is called the {\it canonical extension} of 
$\, M \,$ to $ \P\,$. The invariant $ n $ associated to 
$ M $ can be defined cohomologically in terms of 
its canonical extension, namely 
$$ 
n := \dim_{\C} H^{1}(\P, \ms{M}(-1))\ .
$$
By analogy with the commutative case 
(see \cite{N}, Chapter 2), it is suggestive to think 
of (and refer to) $\, n \,$ as the {\it `second Chern class'} 
$\, c_2(\ms{M}) \,$ of the sheaf $ \ms{M}\,$; however, 
in general, it seems that Chern classes have yet to be defined 
in the realm of noncommutative projective geometry.
\begin{lemma}
\la{L2} 
Let $\, \underline{\mbox{\rm Hom}}\,(\ms{M},\ms{M}) := 
\bigoplus_{k \in \Z} \mbox{\rm Hom}\,(\ms{M},\ms{M}(k)) \,$
be the (graded) endomorphism algebra of the canonical extension
$ \ms{M} $ of $ M $ to $ \P\,$. 
Then, for $ k \gg 0 \,$, we have isomorphisms of 
graded vector spaces: $\, \bs{D}_{\geq k} \cong 
\underline{\mbox{\rm Hom}}(\ms{M},\ms{M})_{\geq k}\,$.
In particular,
\begin{equation}
\la{60}
\dim\, D_{k} = \dim\,\mbox{\rm Hom}(\ms{M},\ms{M}(k)) 
\quad \forall k \gg 0 \ .
\end{equation}
\end{lemma}
\begin{proof}
Without loss of generality, we may identify $ M $ with one of 
the distinguished (fractional) ideals of $ A \,$ (see \cite{BW3},
Section~5), and define its canonical extension $ \ms{M} $ to $ \P $ 
with the help of the induced filtration $\,M_{k} = M \cap Q_{k} \,$.
Then, by (\ref{1}) and (\ref{11}), we have $\, D_{k} = \{\, q \in Q\ :\ 
q M_{i} \subseteq M_{i+k}\ \mbox{for all}\ i\,\}\,$, and hence
\begin{equation}
\la{7}
\bs{D}_{k} \cong 
\mbox{\rm Hom}_{\bs{A}}(\bs{M},\bs{M}(k))\quad \mbox{for all}\ 
k \in Z \ ,
\end{equation}
where $ \mbox{\rm Hom}_{\bs{A}} $ 
in taken in the category $ \Gma\,$.

On the other hand, 
\begin{equation}
\la{8}
\mbox{\rm Hom}\,(\ms{M},\ms{M}(k)) = 
\lim\limits_{\longrightarrow}^{}\,
\mbox{\rm Hom}_{\bs{A}}(\bs{M}_{\geq p},\bs{M}(k))\ .
\end{equation}
To calculate the inductive limit in (\ref{8}) 
we fix some $\, p \gg 0 \,$ and consider the short exact 
sequence in $ \Gma\,$:
\begin{equation}
\la{80}
0 \to \bs{M}_{\geq p} \to \bs{M} \to \bs{M}/\bs{M}_{\geq p} \to 0\ .
\end{equation}
Dualizing (\ref{80}) with $ \bs{M}(k) $ we have
\begin{eqnarray}
0 &\to & 
\mbox{\rm Hom}_{\bs{A}}(\bs{M}/\bs{M}_{\geq p},\bs{M}(k)) \to 
\mbox{\rm Hom}_{\bs{A}}(\bs{M},\bs{M}(k)) \to \nonumber \\*[1ex]
& \to & \mbox{\rm Hom}_{\bs{A}}(\bs{M}_{\geq p},\bs{M}(k)) \to
\mbox{\rm Ext}_{\bs{A}}^{1}(\bs{M}/\bs{M}_{\geq p},\bs{M}(k)) \to 
\ldots \nonumber
\end{eqnarray}
Now, since $ \bs{M}/\bs{M}_{\geq p} $ is torsion
whereas $ \bs{M}(k) $ is torsion-free, 
$\, \mbox{\rm Hom}_{\bs{A}}(\bs{M}/\bs{M}_{\geq p},\bs{M}(k))\,$
must be zero for all $k\,$. On the other hand, $\,
\mbox{\rm Ext}_{\bs{A}}^{1}(\bs{M}/\bs{M}_{\geq p},\bs{M}(k)) = 0 \,$
for $\, k \geq k_{0} \gg 0 \,$ (where $ k_0 $ is independent of
$ p $) because the algebra $ \bs{A} $ satisfies the 
$ \chi$-condition of Artin-Zhang (see \cite{AZ}, Proposition 3.5(1)). 
It follows that
\begin{equation}
\la{9}
\mbox{\rm Hom}_{\bs{A}}(\bs{M}_{\geq p},\bs{M}(k)) \cong
\mbox{\rm Hom}_{\bs{A}}(\bs{M},\bs{M}(k)) \quad \forall k \gg 0\ .
\end{equation}
Combining (\ref{7}) with (\ref{8}) and (\ref{9}) we get the result. 
\end{proof}

\remark\ 
It seems quite plausible that the 
entire vector spaces $ \bs{D} $ and 
$ \underline{\mbox{\rm Hom}}\,(\ms{M},\ms{M}) $ are isomorphic. 
But for our purposes the result of Lemma~\ref{L2} will suffice.

\vspace{1ex}

To finish our calculation it remains to evaluate 
$\, \dim\,\mbox{\rm Hom}\,(\ms{M},\ms{M}(k)) \,$ for
$ k \gg 0\,$ and compare the result with formula 
(\ref{3}). This can be done in many different 
ways of which we will choose (hopefully) the most elementary 
one. 

First, following \cite{KKO} (see {\it loc. cit.}, Section~5.3), 
we introduce  the {\it dual sheaf} $\, \ms{M}^{\vee} := 
\underline{\ms{H}om}\,(\ms{M}, \ms{O})\, $ of $ \ms{M} $
as an object of the quotient category of {\it left} graded 
modules\footnote{For simplicity, 
we will use the same notation $ \P = \mbox{\sf Proj}(\bs{A})\,$, 
$\, \ms{O} = \widetilde{\bs{A}} \,$, etc. for left and right objects.}
over $ \bs{A}\,$. 
\begin{lemma}
\la{L3} 
If $ \ms{M}  \in \coh(\P) $ is trivial over the line
at infinity in $ \P \,$ so is $ \ms{M}^{\vee} \,$.
In that case the twisted sheaves $ \ms{M}(k) $ and 
$ \ms{M}^{\vee}(k) $ have the same Euler characteristics
for all $ k \in \Z\,$:
\begin{equation}
\la{10}
\chi(\P,\, \ms{M}^{\vee}(k)) = \chi(\P,\, \ms{M}(k)) = 
\frac{1}{2}(k+1)(k+2) - n\ .
\end{equation}
\end{lemma}
\begin{proof}
Since $ \ms{M}^{\vee} $ is a torsion-free sheaf of rank $1\,$, 
the first statement of the lemma follows from the second. 
On the other hand, in (\ref{10}) we need only to prove the first 
equality (the second follows immediately from \cite{BW3}, Theorem~4.5).
Furthermore, since 
$\, \chi(\P,\, \ms{M}^{\vee}(k)) \,$ depends polynomially on $ k \,$, 
we may assume $\, k \gg 0 \, $ in which case we have
\begin{equation}
\la{12}
\chi(\P,\, \ms{M}^{\vee}(k)) = 
\dim \,\mbox{\rm Hom}\,(\ms{M},\ms{O}(k)) = 
\dim\, H^{2}(\P,\,\ms{M}(-k-3))\ .
\end{equation}
The first equality in (\ref{12}) follows from the definition
of $ \ms{M}^{\vee} $ and Serre's Vanishing theorem;
the second is a consequence of 
Serre's Duality (cf. \cite{BW3}, Theorem~2.4). 
Both theorems of Serre can be applied
in our situation because, as shown in \cite{KKO},
$\, \ms{M} \in \coh(\P) \,$ if and only if $\, \ms{M}^{\vee} \in 
 \coh(\P) \,$. 
Since we know that both $\, H^{0}(\P,\,\ms{M}(d)) \,$  and  
$ \, H^{1}(\P,\,\ms{M}(d)) \,$ vanish for 
$\, d \ll 0 \,$ (see \cite{BW3}, Lemma~10.2), the last term 
in (\ref{12}) coincides 
with the Euler characteristic $\, \chi(\P,\, \ms{M}(-k-3)) \,$.
The result of the lemma follows now immediately from the obvious 
invariance of $\, \chi(\P,\, \ms{M}(k)) \,$ under the change of 
twisting $\, k \mapsto -k-3\,$. 
\end{proof}

\remark\ 
The result of Lemma~\ref{L3} is specific for our noncommutative 
situation. Indeed, the vanishing of $\, H^{1}(\P,\,\ms{M}(d))\,$ for all 
$ d \ll 0 $ is equivalent to the fact that $ \ms{M} $ is a {\it bundle} 
(in the sense of \cite{KKO}) which is not true for extensions of ideals 
to $ \mathbb{P}^2 $ in the commutative case.

\vspace{1ex}

We proceed now to the final step of our calculation.
\begin{lemma}
\la{L5} 
Let $ \ms{M}_1 $ and $ \ms{M}_2 $ be two rank $1$ bundles on $ \P $ 
trivial over the line at infinity. Then, for all $ k \gg 0 \,$ we have
$$
\dim\,\mbox{\rm Hom}\,(\,\ms{M}_1,\ms{M}_2(k)) = 
\frac{1}{2}(k+1)(k+2) - n_1 - n_2 \ ,
$$
where $ n_1 = c_{2}(\ms{M}_1 ) $ and $ n_2 = c_{2}(\ms{M}_2)\,$.
\end{lemma}
\begin{proof}
We start with a locally free resolution of $ \ms{M}_2 $ in $\, \coh(\P) \,$
(see, for example, \cite {LeB}, Corollary~1):
\begin{equation}
\la{13}
0 \to \bigoplus\limits^{m-1}_{i=1} \ms{O}(-k_{i}) \to
\bigoplus\limits^{m}_{j=1} \ms{O}(-l_{j}) \to \ms{M}_2  \to 0\ ,
\quad  k_{i} \, , \, l_{j} \in \Z_{+}\ .
\end{equation}
Shifting (\ref{13}) by $ k $ and applying the 
functor $\, \mbox{\rm Hom}\,(\,\ms{M}_1,\, \mbox{---}\,) \,$
yields
\begin{eqnarray}
\la{14}
0 &\to & \bigoplus\limits^{m-1}_{i=1} 
\mbox{\rm Hom}\,(\,\ms{M}_1,\ms{O}(k-k_{i})) \to
\bigoplus\limits^{m}_{j=1}
\mbox{\rm Hom}\,(\,\ms{M}_1,\ms{O}(k-l_{j}))  \\
&\to & \mbox{\rm Hom}\,(\,\ms{M}_1,\ms{M}_2(k))
\to \bigoplus\limits^{m-1}_{i=1} 
\mbox{\rm Ext}^{1}\,(\,\ms{M}_1,\ms{O}(k-k_{i})) \to \ldots \ .\nonumber
\end{eqnarray}
If $ k \gg 0 $ the Ext term in (\ref{14}) is zero because
$ \ms{M}_1 $ is a bundle. Therefore, in that case we have
\begin{eqnarray}
\la{15} 
\qquad \lefteqn{\dim\,\mbox{\rm Hom}(\ms{M}_1,\ms{M}_2(k)) =} \\
&& =\sum\limits^{m}_{j=1}\dim\,\mbox{\rm Hom}(\ms{M}_1,\ms{O}(k-l_{j})) 
- \sum\limits^{m-1}_{i=1}\dim\,\mbox{\rm Hom}(\ms{M}_1,\ms{O}(k-k_{i}))\ .
\nonumber
\end{eqnarray}
By Serre's Vanishing theorem we may replace dimensions of 
each Hom term in the right-hand side of (\ref{15}) by the Euler 
characteristics of the corresponding twists of the dual sheaf 
$ \ms{M}^{\vee}_{1}\,$ and then use repeatedly Lemma~\ref{L3}. 
Thus, if  $\, k \gg 0 \,$, we have 
\begin{eqnarray}
\la{16}
\dim\,\mbox{\rm Hom}\,(\,\ms{M}_1,\ms{M}_2(k)) 
&=& 
\sum\limits^{m}_{j=1} \chi(\P,\, \ms{M}^{\vee}_1(k-l_{j})) -
\sum\limits^{m-1}_{i=1} \chi(\P,\, \ms{M}^{\vee}_1(k-k_{i})) 
\nonumber\\*[1ex]
&=& 
\sum\limits^{m}_{j=1} \chi(\P,\, \ms{M}_1(k-l_{j})) -
\sum\limits^{m-1}_{i=1} \chi(\P,\, \ms{M}_1(k-k_{i}))  
\nonumber \\*[1ex]
&=& 
\left(\sum\limits^{m}_{j=1} \chi(\P,\, \ms{O}(k-l_{j})) - mn_1\right) - 
\nonumber \\*[1ex]
&-& \left(\sum\limits^{m-1}_{i=1} \chi(\P,\, \ms{O}(k-k_{i})) - 
(m-1)n_1\right)  \nonumber \\*[1ex]
&=& 
\left(\sum\limits^{m}_{j=1} \chi(\P,\, \ms{O}(k-l_{j}))
-\sum\limits^{m-1}_{i=1} \chi(\P,\, \ms{O}(k-k_{i}))\right) - n_1  
\nonumber \\*[1ex]
&=& \chi(\P,\, \ms{M}_2(k)) - n_1 = \frac{1}{2}(k+1)(k+2) - n_1 - n_2 \ . 
\nonumber
\end{eqnarray}
Note the expression in last parentheses being equal to 
$\, \chi(\P,\, \ms{M}_2(k)) \,$ follows immediately from
the exact sequence (\ref{13}).
\end{proof}

Finally, taking $\, \ms{M}_1 = \ms{M}_2 = \ms{M} $ and 
comparing Lemma~\ref{L5} with formulas (\ref{3}) and (\ref{60}) of 
Lemmas~\ref{L1} and~\ref{L2} we arrive at our main 
\begin{theorem}
\la{T2}
Let $ M $ be a finitely generated rank $1$ torsion-free $A$-module
whose canonical extension to $ \P $ has the `second Chern class'
$\, n \,$. Let $ D := \mbox{\rm End}_{A}(M) \,$.  Then the 
LM-invariant of the algebra $ D $ is given by the formula 
$\, p_{D} = 2 n \,$.
\end{theorem}

\section{Generalizations}

We would like to mention two natural generalizations suggested 
by the results of this paper. First, Lemma~\ref{L5} hints at 
the possibility of defining a {\it relative} version of 
LM-invariants: given two ideal classes in $ A \,$ with distinguished 
representatives $ M_1 $ and $ M_2 \,$ (say), 
we may identify $\, \mbox{\rm Hom}_{A}(M_1, M_2) \,$ 
with a subspace in $ Q $ (as in (\ref{1})), and then 
define the number $\, p_{12} \,$ as the codimension of 
the associated graded space $\, \bs{G}\mbox{\bf Hom}_{A}(M_1, M_2)\,$ in 
$\,\bs{GA} \,$. The argument similar to that of 
\cite{LM}, Proposition~2.4, shows $\, p_{12} < \infty \,$, 
and Lemma~\ref{L5} then immediately implies the relation 
$\,p_{12} = n_1 + n_2 \,$, where $\, n_1 = c_{2}(\ms{M}_1) \,$ 
and $\, n_2 = c_{2}(\ms{M}_2) \,$.

Second, the interpretation of LM-invariants as arithmetic 
genera of noncommutative projective surfaces (Lemma~\ref{L1})
suggests an obvious generalization to higher dimensions.
As in the original situation of \cite{LM},
one would expect these generalized LM-invariants 
to be useful in the study of rings of differential operators 
on singular rational algebraic varieties. 
An interesting problem in this direction 
would be to compute the LM-type invariants for the algebras of 
differential operators on varieties of quasi-invariants 
of finite reflection groups (see \cite{BEG}).

\bibliographystyle{amsalpha}

\end{document}